\documentclass[pdflatex,sn-mathphys-num]{sn-jnl}

\usepackage{graphicx}%
\usepackage{multirow}%
\usepackage{amsmath,amssymb,amsfonts}%
\usepackage{amsthm}%
\usepackage{mathrsfs}%
\usepackage[title]{appendix}%
\usepackage{xcolor}%
\usepackage{textcomp}%
\usepackage{manyfoot}%
\usepackage{booktabs}%
\usepackage{algorithm}%
\usepackage{algorithmicx}%
\usepackage{algpseudocode}%
\usepackage{listings}%

\theoremstyle{thmstyleone}%
\newtheorem{theorem}{Theorem}
\newtheorem{remark}[theorem]{Remark}
\newtheorem{proposition}[theorem]{Proposition}%
\newtheorem{corollary}[theorem]{Corollary}

\theoremstyle{thmstyletwo}%

\theoremstyle{thmstylethree}%

\begin{document}

\title[The $n$-th Zb\'{a}ganu Constant]{The $n$-th Zb\'{a}ganu Constant}


\author[1,2]{\fnm{Ifronika} \sur{}}\email{ifronika@itb.ac.id}
\equalcont{These authors contributed equally to this work.}

\author*[2]{\fnm{Hendra} \sur{Gunawan}}\email{hgunawan@itb.ac.id}
\equalcont{These authors contributed equally to this work.}

\author[2]{\fnm{Janny} \sur{Lindiarni}}\email{janny@itb.ac.id}
\equalcont{These authors contributed equally to this work.}

\affil[1]{\orgdiv{Doctoral Program of Mathematics}, \orgname{Faculty of Mathematics and Natural Sciences, Bandung Institute of Technology}, \orgaddress{\street{Jl. Ganesha 10}, \city{Bandung}, \postcode{40132}, \country{Indonesia}}}

\affil*[2]{\orgdiv{Department of Mathematics}, \orgname{Faculty of Mathematics and Natural Sciences, Bandung Institute of Technology}, \orgaddress{\street{Jl. Ganesha 10}, \city{Bandung}, \postcode{40132}, \country{Indonesia}}}


\abstract{In this paper, we introduce the $n$-th Zb\'{a}ganu constant for Banach spaces 
and investigate its basic properties. We establish some inequalities which give the 
relationship between the $n$-th Zb\'{a}ganu constant and two other constants, namely 
the $n$-th Von Neumann-Jordan constant and the $n$-th James constant. In addition, we 
determine the values of the $n$-th Zb\'{a}ganu constant for Lebesgue spaces, Morrey 
spaces, and discrete Morrey spaces.}




\maketitle

\section{Introduction}\label{sec1}

One way to study the convexity of a Banach space is through its geometric constants,
such as the Von Neumann-Jordan constant or the James constant. For a Banach space 
$(X,\|\cdot\|_X)$, the Von Neumann-Jordan constant for $X$ is defined in \cite{Jordan} by
\begin{equation*}
C_{NJ}(X):=\sup\bigg\{
\frac{\|x+y\|_X^2+\|x-y\|_X^2}{2(\|x\|_X^2+\|y\|_X^2)} : x,y\in X \backslash \{ 0 \} 
\bigg\}.
\end{equation*}
In general, we have $1 \le C_{NJ}(X) \le 2$, where the minimum value, that is $C_{NJ}(X)
=1$, is attained if and only if $X$ is a Hilbert space. The value of the Von 
Neumann-Jordan constant $C_{NJ}(X)$ is closely related to the convexity of $X$, as 
investigated by Kato and Takahashi \cite{Kato}. Further results on this 
constant and its relationship with other geometric properties can be found in 
\cite{Takahashi,Alonso,Gunawan5,Cui,Yang}.

In 1998, Kato \textit{et al.} \cite{Kato2} introduced a higher-order version of the Von 
Neumann-Jordan constant. For $n \ge 2$, the $n$-th Von Neumann-Jordan constant 
$C_{NJ}^{(n)}(X)$ is defined by
\begin{equation*}
C_{NJ}^{(n)}(X):=\sup\left\{
\frac{\sum\limits_{\theta_i=\pm 1} \left\|x_1+\sum\limits_{i=2}^n \theta_ix_i 
\right\|_X^2}{2^{n-1}\sum\limits_{i=1}^n \|x_i\|_X^2} : x_i\in X, \sum\limits_{i=1}^n 
\|x_i\|_X^2 \neq 0, i=1,\ldots, n \right\}.
\end{equation*}
Kato \textit{et al.} then showed that, for every Banach space $X$, the inequality
$1 \le C_{NJ}^{(n)}(X)\le n$ holds. In particular, $C_{NJ}^{(n)}(X)=1$ if and only if
$X$ is a Hilbert space \cite{Kato2}.


Another important geometric constant is the James constant, which was developed based
on James's work in \cite{James}. For a Banach space 
$(X,\|\cdot\|_X)$, let $S_X := \{x \in X : \|x\|_X = 1\}$, the unit sphere in $X$. 
The James constant is defined in \cite{Gao} by
\begin{equation*}
C_J(X) := \sup \{\min(||x+y||_X, ||x-y||_X) : x,y \in S_X\}.
\end{equation*}
A higher-order version of the James constant was introduced by Maligranda \textit{et al.} \cite{Maligranda}. For $n\ge 2$, the $n$-th James constant is defined by
\begin{equation*}
C_J^{(n)}(X) := \sup \left\{\min\limits_{\theta_i=\pm1} \left\{\left\|x_1+ 
\sum\limits_{i=2}^n \theta_ix_i \right\|_{X} \right\} : x_i \in S_X, 
i=1,\ldots,n\right\}.
\end{equation*}
If $\dim(X)\ge n$, then the inequality $\sqrt{n}\le C_J^{(n)}(X) \le n$ holds.
Moreover, if $X$ is a Hilbert space, then $C_J^{(n)}(X)=\sqrt{n}$. The converse,
however, is false \cite{Maligranda}.

In 2002, Zb\'{a}ganu \cite{Zbaganu,Mizuguchi} introduced a different constant, namely
\begin{equation*}\label{eq:4} 
C_{Z}(X) :=\sup\bigg\{
\frac{\|x+y\|_X\|x-y\|_X}{\|x\|_X^2 +\|y\|_X^2} : x,y\in X\setminus\{0\} \bigg\}.
\end{equation*}
One may observe that $1\le \frac12 C_J^2(X)\le C_Z(X)\leq C_{NJ}(X) \leq 2.$
Note that if $X$ is a Hilbert space, then $\frac12 C_J^2(X)=C_Z(X)=C_{NJ}(X)=1$.
A modified version of the Zb\'{a}ganu constant is studied in \cite{Li}. Some
applications of the Zb\'{a}ganu constant may be found in \cite{Dinarvand}.

Motivated by the higher-order Von Neumann-Jordan and James constants, we introduce the 
$n$-th Zb\'{a}ganu constant and study its basic properties. [We are surprised that
no one has studied this constant after Zb\'{a}ganu's work in 2002.]
After defining the constant, we present its relationship 
with the $n$-th Von Neumann-Jordan constant and the $n$-th James constant. 
In addition, we determine the value of the $n$-th Zb\'{a}ganu constant for some
classical spaces, including Morrey spaces and discrete Morrey spaces. We find that
the value of the $n$-th Zb\'{a}ganu constant for each of those spaces is equal to 
the values of the $n$-th Von Neumann-Jordan constant and the $n$-th James constant.

\section{The $n$-th Zb\'{a}ganu Constant}\label{sec2}


Let $X$ be a Banach space. For $n \ge 2$, we define the $n$-th Zb\'{a}ganu constant by
\[
C_Z^{(n)}(X)= \sup  \left\{ \frac{\prod\limits_{\theta_i=\pm 1} \left\|x_1+ 
\sum\limits_{i=2}^n \theta_ix_i \right\|_X^{\frac{1}{2^{n-2}}} }{\sum\limits_{i=1}^n 
\|x_i\|_X^2} : x_i \in X\ \text{for}\ i=1, \ldots,n,\ \sum\limits_{i=1}^n \|x_i\|_X^2 
\neq 0\right\}.
\]
Note that, for $n=2$, the formula reduces to that for the original Zb\'{a}ganu constant, 
that is, $C_Z^{(2)}(X)=C_Z(X)$.

One may observe that $C_Z^{(n)}(X) \ge 1$ for every $n\ge 2$. To see this, let $x_1, 
\ldots, x_n \in X$. Choose $x_1\neq0$ and $x_i=0$ for $i=2,\ldots n.$ Since there are 
$2^{n-1}$ choices of signs, we have 
\[
\prod\limits_{\theta_i=\pm 1} \left\|x_1+ \sum\limits_{i=2}^n \theta_ix_i\right\|_X=
\|x_1\|_X^{2^{n-1}}.
\]
Therefore, 
\begin{align*}
\frac{\prod\limits_{\theta_i=\pm 1} \left\|x_1+ \sum\limits_{i=2}^n \theta_ix_i 
\right\|_X^{\frac{1}{2^{n-2}}} }{\sum\limits_{i=1}^n \|x_i\|_X^2} =
\frac{\left(\|x_1\|_X^{2^{n-1}}\right)^{\frac{1}{2^{n-2}}}}{\|x_1\|_X^2}=1.
\end{align*}
By the definition of $C_Z^{(n)}(X)$, we obtain $C_Z^{(n)}(X) \ge 1$, as claimed. 

\medskip

The following theorem gives the relationship between the $n$-th Zb\'{a}ganu constant,
the $n$-th Von Neumann-Jordan constant, and the $n$-th James constant.

\begin{theorem}\label{th:1}
For any Banach space $X$, we have 
\[[C_J^{(n)}(X)]^2 \le nC_{Z}^{(n)}(X) \le nC_{NJ}^{(n)}(X),\] 
for every $n \ge 2$.
\end{theorem}

\begin{proof}
We shall first show that $[C_J^{(n)}(X)]^2 \le nC_{Z}^{(n)}(X)$ for every $n\ge2$.
Let $x_1,\ldots,x_n \in S_X$ and $m=\min\limits_{\theta_i=\pm1} \left\{\left\|x_1+ 
\sum\limits_{i=2}^n \theta_ix_i \right\|_{X} \right\}$. Then,
\begin{align*}
    m &\le \left( \prod\limits_{\theta=\pm1}\left\|x_i+ \sum\limits_{i=2}^n \theta_ix_i 
    \right\|_{X}\right)^{\frac{1}{2^{n-1}}}\\
    &=\left(\frac{n\prod\limits_{\theta=\pm1}\left\|x_i+ \sum\limits_{i=2}^n 
    \theta_ix_i\right\|^{\frac{1}{2^{n-2}}}_{X}}{\sum\limits_{i=1}^n 
    \|x_i\|_X^2}\right)^{\frac{2^{n-2}}{2^{n-1}}}\\
    &=[nC_Z^{(n)}(X)]^\frac12.
\end{align*}
Taking the supremum over all $x_1,\ldots,x_n \in S_X,$ and then squaring both sides,
we obtain 
\[[C_J^{(n)}(X)]^2\le n C_Z^{(n)}(X).\]

We shall now establish the inequality $C_{Z}^{(n)}(X) \le C_{NJ}^{(n)}(X)$ 
for every $n \ge 2$.
Let $x_1, \ldots, x_n \in X$ satisfy $\displaystyle \sum\limits_{i=1}^n \|x_i\|_X^2 >0$. 
By virtue of the GM-QM Inequality we have, 
\[
   \left( \prod_{\theta_i = \pm 1} \left\|x_1+ \sum\limits_{i=2}^n \theta_i x_i 
   \right\|_X \right)^{\frac{1}{2^{n-1}}} \le \sqrt{\frac{\sum\limits_{\theta_i=\pm 1} 
   \left\|x_1+ \sum\limits_{i=2}^n \theta_i x_i  \right\|_X^2}{2^{n-1}}},
\]
whence
\[
  \left( \prod_{\theta_i = \pm 1} \left\|x_1+ \sum\limits_{i=2}^n \theta_i x_i 
  \right\|_X \right)^{\frac{1}{2^{n-2}}} \le \frac{\sum\limits_{\theta._i=\pm 1} 
  \left\|x_1+ \sum\limits_{i=2}^n \theta_i x_i  \right\|_X^2}{2^{n-1}}.  
\]
Dividing both sides by $\displaystyle \sum_{i=1}^n \|x_i\|_X^2 >0,$ we get
\[
\frac{\left( \prod\limits_{\theta_i = \pm 1} \left\| x_1+\sum\limits_{i=2}^n 
\theta_i x_i \right\|_X \right)^{\frac{1}{2^{n-2}}}}{\sum\limits_{i=1}^n \|x_i\|_X^2} 
\le \frac{\sum\limits_{\theta_i=\pm 1} \left\| x_1+\sum\limits_{i=2}^n \theta_i x_i  
\right\|_X^2}{2^{n-1}\sum\limits_{i=1}^n \|x_i\|_X^2} \le C_{NJ}^{(n)}(X).
\]
By taking the supremum over all $x_1, \ldots, x_n \in X$, we obtain
\[
C^{(n)}_Z(X) \le C_{NJ}^{(n)}(X),
\]
as desired.

Combining this with the earlier inequality, the proof is complete.
\end{proof}

\begin{corollary}
For every $n\ge 2$, we have $1\le C_Z^{(n)}(X)\le n$. Moreover, if $X$ is a Hilbert
space, then $C_Z^{(n)}(X)=1$.
\end{corollary}

\section{Some Examples}

As our first example, we determine the value of the $n$-th Zb\'{a}ganu constant for
the classical $L^p$ and $\ell^p$ spaces, where $1\le p\le 2$.

\begin{theorem}
For $X=L^p$ and $X=\ell^p$, where $1\le p\le 2$, we have $C_Z^{(n)}(X)=n^{2/p}-1$.
\end{theorem}
\begin{proof}
    For $X=L^p$ and $X=\ell^p$, where $1\le p\le 2$, we know that
    $C_{NJ}^{(n)}(X)=n^{2/p-1}$ \cite{Kato} and $C_J^{(n)}(X)=n^{1/p}$ \cite{Maligranda}.
    It thus follows from Theorem \ref{th:1} that $C_Z^{(n)}(X)=n^{2/p}-1$.
\end{proof}

\medskip

The following proposition provides a sufficient condition for which each of the three 
constants attains its maximum value, which is $n$. The proof follows immediately from 
the definitions and the fact that each of the three constants is bounded by $n$.

\begin{proposition}\label{prop:2} 
Let $X$ be a Banach space. If there exist $x_1,\ldots x_n \in X$ with 
$\|x_1\|_{X}=\cdots=\|x_n\|_{X}=1$ such that, for every choice of $\theta_i=\pm1$ 
$(i=2,\ldots,n)$,
\[
\|x_1+ \sum\limits_{i=2}^n \theta_ix_i\|_{X}= n,
\]
then $C_J^{(n)}(X)=C_Z^{(n)}(X)=C_{NJ}^{(n)}(X)=n$.
\end{proposition}

The above proposition has been employed to obtain that the $n$-th Von Neumann-Jordan 
constant and the $n$-th James constant for Morrey spaces 
$\mathcal{M}^p_q=\mathcal{M}^p_q(\mathbb{R}^d)$ and 
discrete Morrey spaces $\ell^p_q=\ell^p_q(\mathbb{Z}^d)$, where $1\le p<q<\infty$,
are both equal to $n$ \cite{Gunawan5}. The same construction of functions and sequences 
in $\mathcal{M}^p_q$ and $\ell^p_q$, respectively, as in \cite{Gunawan5}, gives us the 
following result for the $n$-th Zb\'{a}ganu constant, as our next example.

\begin{theorem}
For $1\le p<q<\infty$ and $n\ge 2$, we have 
$$C_{Z}^{(n)}(\mathcal{M}^p_q)=C_{Z}^{(n)}(\ell^p_q)=n.$$
\end{theorem}

\begin{remark}
    The above result also follows from Theorem \ref{th:1} since we already know that
    $C_J^{(n)}(X)=C_{NJ}^{(n)}(X)=n$ for $X=\mathcal{M}^p_q$ and $X=\ell^p_q$.
\end{remark}

\section*{Acknowledgment}
The research is supported by P2MI-ITB 2026 Program.

\section*{Declarations}
The authors declare that they do not have conflict of interest.

\end{document}